\documentclass[reqno,12pt]{amsart}
\usepackage{amsmath,amssymb,tabularx,setspace,color}
\usepackage{multirow}
\newtheorem{theorem}{Theorem}[section]

\newtheorem{definition}{Definition}[section]
\newtheorem{remark}{Remark}[section]
\newtheorem{corollary}{Corollary}[section]
\makeatletter
\renewcommand\section{\@startsection {section}{1}{\z@}%
                                   {-3.5ex \@plus -1ex \@minus -.2ex}%
                                   {2.3ex \@plus.2ex}%
                                   {\normalfont\large\bfseries}}
\makeatother

\usepackage{fancyhdr}

\begin{document}
\doublespace

\title{An exact test for  renewal increasing mean residual life}

\author{S\lowercase{udheesh} K K\lowercase{attumannil}
\\I\lowercase{ndian} S\lowercase{tatistical} I\lowercase{nstitute},
  C\lowercase{hennai}, I\lowercase{ndia}}

\maketitle

\begin{abstract}
In this paper, we develop an exact test for testing exponentiality against renewal increasing mean residual life class. Pitman's asymptotic efficacy value shows that our test perform well. Some numerical results are presented to demonstrate the performance of the testing method. We also discuss how the proposed method incorporates the right censored observations.
\\
Keywords: Exponential distribution; Renewal increasing mean residual life; Replacement model; Shock Model; U-statistics.
\end{abstract}

\section{Introduction}



When a device is experiencing random number of shocks
governed by a homogeneous Poisson process, the concept of renewal
increasing mean residual life is very much useful to study the properties of age replacement model.  In this context, test for exponentiality against renewal increasing mean residual life class is used to
determine whether to adopt a planned replace model over unscheduled one. Sepehrifar et al. (2015) developed a non-parametric test against RIMRL$_{shock}$ class and obtained a critical region based on the asymptotic theory of U-statistics. We noted that the critical region developed by Sepehrifar et al. (2015) is incorrect (see Remark 2.1). Motivated by Sepehrifar et al. (2015), we develop an exact test for testing exponentiality against RIMRL$_{shock}$  class. We also obtain the correct critical region of the asymptotic test proposed by Sepehrifar et al. (2015). The case with censored observations also addressed.

The rest of the paper is organized as follows. In Section 2, we propose an exact test for testing exponentiality against RIMRL$_{shock}$ class and then calculate the critical values for different sample sizes. The asymptotic normality proposed test statistic is proved in Section 3. The Pitman's asymptotic efficacy value is also given in this section. In Section 4, we report the result of simulation study carried out to assess the performance of the proposed test. In Section 5, we discuss how to incorporate right censored observation in our study.

\section{Exact Test}
Let $X$ be the lifetime of a device which has absolutely continuous distribution function $F(.)$. Suppose $\bar{F}(x)=P(X>x)$ denotes the survival function of $X$ at $x$. Also let $\mu=E(X)=\int_{0}^{\infty}\bar{F}(t)dt<\infty$. Assume that the device under consideration is experiencing a random shock. Suppose $N(t)$ denotes the total number of shocks
up to time $t$ with probability mass function $P(N(t)=j)=F^{j}(t)-F^{j+1}(t),j=0,1,2,....$
Suppose that the random variable $W_j, j=0,1,2,… $ quantify the amount of hidden lifetime absorbed by the $j$th shock with $W_0=0$ and having common distribution function $G(x)=P(W_j\leq x)$.  The total cumulative life damage up to time $t$ is defined as $Z(t)=\sum_{j=0}^{N(t)}W_j$ with the cumulative distribution function $Q(x)=P(Z(t)\leq x)=\sum_{j=0}^{\infty}G^{(j)}[F^{(j)}(t)-F^{(j+1)}(t)]$. It is assumed that the unit fails when the total life-damage exceeds a pre-specified level $x>0$.
We refer to  Glynn and Whitt (1993), Roginsky (1994) and Sepehrifar et al. (2015) for discussion related this framework.

 Let $X^*=X-Z(t)$ be the residual lifetime of an operating device with cumulative damage $Z(t)$. Note that the realizations of $X^*$ is available to us for further analysis. Consider a device subjected to $N(t)$ number of shocks up to time $t$. Given that such a device is in an operating situation at time instant $t$ after installation, the MRL function of $X^*$ denoted by $m^*(t)$ is defined by $m^*(t)=E(X^*-t|X^*\geq t)$. Note that the total life-damage will not exceed the threshold level $x$. From the definitions it is evident that the random variables $X^*$ and $Z(t)$ are independent. Next we give the definition of RIMRL$_{shock}$ class (Sepehrifar et al., 2015).
\begin{definition}
The mean residual life of a device under shock model (MRL$_{shock}$) at time $t$ is defined as $$m^*(t)=\frac{1}{\bar{r}(t)}\int_t^{\infty}\bar{r}(z)dz,$$ where $\bar{r}(z)=\int_0^x\bar{F}(z+w)dQ(w)$.
\end{definition}
\begin{definition}
The random variable $X$ belongs to the RIMRL$_{shock}$ class if the function $m^*(t)$ is a non-decreasing function for all $t>0$.
\end{definition}
 We are interested to test the null hypothesis
 \begin{equation*}
  H_{0}: F^* \quad \text{is exponential}
\end{equation*}
against the alternatives
\begin{equation*}
  H_{1} : F^* \quad \text{is RIMRL$_{shock}$ (and not exponential)},
\end{equation*}on the basis of a random sample $X_1^*,X_2^*,...,X_n^*$; from an absolutely continuous distribution function $F^*$. For the above testing problem Sepehrifar et al. (2015) proposed a non-parametric test based on the departure measure $\Delta^*(F^*)$ defined by
\begin{equation*}
  \Delta^*(F^*)=\frac{1}{\mu^*}E_{f^*}(min(X_1^*,X_2^*)-\frac{1}{2}X_1^*)=\frac{\Delta(F^*)}{\mu^*},
\end{equation*}where $\Delta(F^*)=E_{f^*}(min(X_1^*,X_2^*)-\frac{1}{2}X_1^*)$  and $\mu^*=E(X_{1}^*)$. Based on U-statistics theory Sepehrifar et al. (2015) obtained the following test statistic
%
%
\begin{equation}\label{eq2.1}
  \widehat{\Delta}^*=\frac{\widehat{\Delta}}{\bar{X}^*},
\end{equation}where $\bar{X}^*=\frac{1}{n}\sum_{i=1}^{n}X^*_i$ and $\widehat{\Delta}=\frac{2}{n(n-1)}\sum_{i=1}^{n}\sum_{j<i;j=1}^{n}h(X^{*}_{i},X^{*}_{j})$ with $h(X_1^*,X_2^*)=min(X_1^*,X_2^*)-\frac{1}{2}X_1^*$. Hence the test procedure is to reject the null hypothesis $H_0$ in favour of $H_1$ for large values of $\widehat{\Delta}^*$.

\begin{remark}
For the testing problem discussed here, Sepehrifar et al. (2015) obtained a critical region based on the asymptotic variance $\frac{7}{48}$. However as we shown (see Section 3)  the asymptotic variance is $\frac{1}{12}$. We could not find enough details in their paper to explain the discrepancy.
\end{remark}


Motivated by this discrepancy, next we develop an exact test based on the test statistics $\widehat{\Delta}^*$ and calculate the critical values for different sample size. We use a result due to Box (1954) to find the exact null distribution of the test statistic.

%
\begin{theorem}Let $X^*$ be continuous non-negative random variable with $\bar{F}^*(x)=e^{- \frac{x}{2}}$. Let $X_{1}^*,X_{2}^*, ...,X_{n}^*$ be independent and identical samples from $F^*$. Then for fixed $n$
\begin{equation*}
  P(\widehat{\Delta}^*> x)=\sum_{i=1}^{n}\prod_{j=1,j\neq i}^{n}\Big(\frac{d_{i,n}-x}{d_{i,n}-d_{j,n}}\Big)I(x,d_{i,n}),
\end{equation*}provided $d_{i,n}\neq d_{j,n}$ for $i\neq j$, where
\[I(x,y) = \left\{ {\begin{array}{*{20}{c}}
{1{\rm{ }}\quad if{\rm{ }}\quad x \leq y}\\
{ 0{\rm{ }}\quad if{\rm{ }}\quad x > y}
\end{array}} \quad \textit{and}\quad  d_{i,n}= \frac{(n-2i+1)}{2(n-1)}.\right. \]
\end{theorem}
 \noindent {\bf{Proof:}}
  First we express the test statistics in terms of order statistics.
   Note that
   \begin{equation*}
     \frac{2}{n(n-1)}\sum_{i=1}^{n}\sum_{j<i;j=1}^{n}min(X^{*}_{i},X^{*}_{j})=\frac{2}{n(n-1)}\sum_{i=1}^{n}(n-i)X_{i},
   \end{equation*}
  where $X^*_{(i)}$, $i=1,2,...,n$, is the $i$-th order statistics based on the random sample $X^*_{1}, X^*_{2},...X^*_{n}$; from $F^*$.  After some algebraic manipulation, we can express the equation (\ref{eq2.1}) as
\begin{eqnarray}\label{eq2.10}
  \widehat{\Delta}^*&=&\frac{ \frac{1}{2n(n-1)}\sum_{i=1}^{n}(3n-4i+1)X^*_{(i)}}{\bar{X}^*}.
\end{eqnarray}
 Rewrite the denominator of the equation (\ref{eq2.10}) as
 \begin{equation*}
   \widehat{\Delta}=\sum_{i=1}^{n}X^{*}_{(i)}\Big[\frac{(n-i+1)^{2}}{n(n-1)}-\frac{(n-i)^{2}}{n(n-1)}-\frac{(n+1)}{2n(n-1)}\Big].
 \end{equation*}
 Or \begin{equation*}
   \widehat{\Delta}=\frac{n}{(n-1)}\sum_{i=1}^{n}X^{*}_{(i)}\Big[\frac{(n-i+1)^{2}}{n^{2}}-
   \frac{(n-i)^{2}}{n^{2}}-\frac{(n+1)}{2n^{2}}\Big].
 \end{equation*}
Hence, in terms of the normalized spacings, $D_{i}=(n-i+1)(X_{(i)}^*-X_{(i-1)}^*)$, with $X_{0}^*=0$, we can express the test statistics as
  \begin{equation*}
   \widehat{\Delta}^*=\frac{\sum_{i=1}^{n}d_{i,n}D_{i}}{\sum_{i=1}^{n}D_{i}},
 \end{equation*}
 where $ d_{i,n}$'s  are given by
\begin{eqnarray*}
   d_{i,n}&=&\frac{1}{(n-1)}\big[(n-i+1)-\frac{(n+1)}{2}\big]\\
   &=&  \frac{(n-2i+1)}{2(n-1)}.
\end{eqnarray*}Note that the exponential random variable with rate $\frac{1}{2}$ is distritbuted same as the  $\chi^2$ random variable with 2 degrees of freedom.  Hence the result follows from Theorem 2.4 of Box (1954).

\noindent The critical values of the exact test for different $n$ are tabulated in Table 1.
\begin{center}
\begin{table}[h]\footnotesize
\caption{Critical values of the exact test }
\begin{tabular}{|c|c|c|c|c|}\hline
$n$ &  90\% level &   95\% level &  97.5\% level&  99\% level \\ \hline
2&	0.4000& 	0.4500&	0.4750&	0.4900\\ \hline
3&	0.2764&	0.3419&	0.3883&	0.4292\\ \hline
4&	0.2189&	0.2678&	0.323&	0.3693\\ \hline
5&	0.1883&	0.2383&	0.28&	0.325\\ \hline
6&	0.1679&	0.2131&	0.2508&	0.2927\\ \hline
7&	0.1529&	0.1944&	0.2293&	0.2682\\ \hline
8&	0.1413&	0.1799&	0.2125&	0.2492\\ \hline
9&	0.1319&	0.1682&	0.1989&	0.2336\\ \hline
10&	0.1243&	0.1586&	0.1877&	0.2208\\ \hline
15&	0.0993&	0.1271&	0.1508&	0.178\\ \hline
20&	0.0852&	0.109&	0.1295&	0.1531\\ \hline
25&	0.0758&	0.097&	0.1153&	0.1363\\ \hline
30&	0.0689&	0.0882&	0.1049&	0.1241\\ \hline
40&	0.0594&	0.0761&	0.0905&	0.1072\\ \hline
50&	0.0529&	0.0679&	0.0808&	0.0957\\ \hline
75&	0.0431&	0.0552&	0.0658&	0.078\\ \hline
100& 0.0373& 0.0477& 0.0569& 0.0675\\ \hline
\end{tabular}
\end{table}
\end{center}

 \section{Asymptotic properties}
 In this section, we prove asymptotic normality of the proposed test statistic.  Making use of the asymptotic distribution we also calculate the Pitman's asymptotic efficacy of the test. As mentioned  Sepehrifar et al. (2015) showed that the test statistic has limiting normal distribution, however they  incorrectly stated the asymptotic variance. Hence we give the following results to correct the error occurred in their study.
\begin{theorem} The  distribution of $\sqrt{n}(\widehat{\Delta}-\Delta{(F^*)})$, as $ n \rightarrow \infty $, is Gaussian with mean zero and variance $4\sigma_{1}^{2}$, where $\sigma_1^2$ is the asymptotic variance of $\widehat{\Delta}$ and is given by
\begin{equation}\label{eq3.12}
\sigma_{1}^{2}=\frac{1}{4}Var\Big(2X^*\bar{F^*}(X^*)+2\int_{0}^{X^*}ydF^*(y)-\frac{1}{2}X^*\Big).
\end{equation}
\end{theorem}
\begin{corollary}Let $X^*$ be continuous non-negative random variable with $\bar{F^*}(x)=e^{-\frac{x}{\lambda} }$, then the distribution of $\sqrt{n}\widehat{\Delta}$, as $n\rightarrow\infty$, is Gaussian with mean zero and variance $\sigma_{0}^{2}=\frac{\lambda^{2}}{12}.$
\end{corollary}
%
\begin{corollary}Let $X^*$ be continuous non-negative random variable with $\bar{F^*}(x)=e^{-\frac{x}{\lambda} }$, then the distribution of $\sqrt{n}\widehat{\Delta}^*$, as $n\rightarrow\infty$, is Gaussian with mean zero and variance $\sigma_{0}^{2}=\frac{1}{12}.$
\end{corollary}

\noindent Apart from the exact test we can construct an asymptotic test based on the asymptotic distribution
of  $\widehat{\Delta}^*$. Hence in case of the asymptotic test, for large values of $n$, we reject the null hypothesis $H_{0}$ in favour of the alternative hypothesis $H_{1}$, if
\begin{equation*}
  \sqrt{12n}(\widehat{\Delta}^*)>Z_{\alpha},
  \end{equation*}
where $Z_{\alpha}$ is the upper $\alpha$-percentile of $N(0, 1)$. In fact, this is the correct critical region of the test proposed by Sepehrifar et al. (2015).

Next we study the asymptotic efficiency of the test. The Pitman's asymptotic efficacy is the most frequently used index to make a quantitative comparison of two distinct asymptotic tests for a certain statistical hypothesis.  The Pitman's asymptotic efficacy (PAE) is defined as
\begin{equation*}\label{eq6a}
  PAE(\widehat{\Delta}^*)= \frac{ |\frac{d}{d\lambda}\Delta^*(F^*)|_{\lambda\rightarrow \lambda_{0}}}{ \sigma_{0}},
\end{equation*}where $\lambda_{0}$ is the value of $\lambda$ under $H_{0}$ and $\sigma_{0}^2$ is the asymptotic variance of $\widehat{\Delta}^*$ under $H_{0}$.
In our case, the PAE is given by
\begin{eqnarray*}
  PAE(\Delta^*(F^*))&=& \frac{ |\frac{d}{d\lambda}\Delta^*(F^*)|_{\lambda\rightarrow \lambda_{0}}}{\sigma_{0}} \\
  &=& \sqrt{12}(W'(\lambda_{0})-W(\lambda_{0}) \mu^{*'}_{a}(\lambda_{0})),
\end{eqnarray*}where $W=E(min(X_1^*,X_2^*))$ and $\mu^{*}_{a}$ is the mean of $X^*$ under the alternative hypothesis and the prime denotes the differentiation with respect to $\lambda$.
We calculate the PAE value for three commonly used alternatives which are the members of RIMRL$_{shock}$ class \\
(i) the Weibull family: $\bar{F}^*(x)=e^{-x^{\lambda}}$ for $\lambda>1$, $x\geq 0$\\
(ii) the linear failure rate family: $\bar{F}^*(x)=e^{(-x-\frac{\lambda}{2}x^{2})}$ for $\lambda>0$, $x\geq 0$\\
(iii) the Makeham family: $\bar{F}^*(x)=e^{-x-\lambda(e^{-x}+x-1)}$ for $\lambda>0$, $x\geq 0$.

By direct calculations, we observe that the PAE for Weibull distribution is equal to 1.2005; while for linear
failure rate distribution and the Makeham distribution
these values are, 0.8660 and 0.2828, respectively. 

Next we compare the performance of the proposed test with some other tests available in the context of age replace model by evaluating the PAE of the respective tests. We compare our test with that tests proposed by Li and Xu (2008) and Kayid et al. (2013).  The Table 2 gives the PAE values for different test procedures. From the Table 2, it is clear that our test is quite efficient for the Weibull and linear failure rate alternatives.  Note that the test proposed by Kayid et al. (2013) has good efficacy for Makeham alternative even though their test shows poor performance against the other two given alternatives.

\begin{table}[h]\footnotesize
\caption{Pitman's asymptotic efficacy (PAE)}
\begin{tabular}{|r|r|r|r|r|r|r|}
\hline
      \!Distribution\! &     \!Proposed test\! &     \!Li and Xu (2008)\! & \!Kayid  et al. (2013)\!    \\
\hline
         Weibull\! & 1.2005\! &1.1215\! & 0.4822\! \\
\hline
        Linear failure rate\! &0.8660\!  & 0.5032\! & 0.4564\!   \\
\hline
         Makeham\! & 0.2828\! & 0.2414\! & 2.084\!  \\
\hline
\end{tabular}
\end{table}

\section{Simulation study}
Here we report a simulation study for evaluating the performance of our asymptotic test against various alternatives. The simulation was done using R program.

First we find the empirical type 1 error of the proposed test. Since the
test is scale invariant, we simulate random sample from standard exponential distribution. The simulation is repeated for ten thousand times with different values of $n$ and is reported in Table 3. From the Table 3 it evident that the empirical type 1 error is a very good estimator of the size of the test even for small sample size.

For finding empirical power against various alternatives, we simulate observations from Weibull, linear failure rate and Makeham distributions with different values of $\lambda$ where the distribution functions were given in the Section 3. The empirical powers for the above mentioned alternatives are given in Tables  4, 5 and 6. From these tables we can see that empirical powers of the test approaches one when the $\theta$ values are going away from the null hypothesis value as well as when $n$ takes large values.
\begin{table}
\caption{Empirical type 1 error of the test }
\begin{tabular}{|c|c|c|}\hline
\!$n$\! & \!Type 1  Error (5\% level)\! & \!Type 1  Error (1\% level)\!\\ \hline
10& 0.0635&	0.0123\\ \hline
20& 0.0540&	0.0115\\ \hline
30&	0.0518&	0.0107\\ \hline
40&	0.0520&	0.0110\\ \hline
50&	0.0517&	0.0107\\ \hline
60&	0.0516&	0.0102\\ \hline
70&	0.0515&	0.0102\\ \hline
80&	0.0511&	0.0100\\ \hline
90&	0.0504&	0.0103\\ \hline
100&0.0504&	0.0104\\ \hline
\end{tabular}
\end{table}

\begin{table}[h]
\caption{Empirical Power: Weibull distribution }
\centering
\begin{tabular}{|l|l|l|l|l|l|l|l|l|}
\hline
\multirow{2}{*}{} & \multicolumn{2}{l|}{$\lambda=1.2$} & \multicolumn{2}{l|}{$\lambda=1.4$} & \multicolumn{2}{l|}{$\lambda=1.6$} & \multicolumn{2}{l|}{$\lambda=1.8$} \\ \cline{2-9}
 $n$& 5\!\% & 1\%& 5\!\% & 1\%& 5\% & 1\!\%&  5\!\%& 1\%   \\ \hline
 60\!&	0.50&	0.23\!&	0.93&	0.76&	0.99\!&	0.97&	1.00&	0.99\!\\ \hline
 70\!&    0.55&	0.27\!&	0.96&	0.84&	0.99\!&	0.99&	1.00&	1.00\!\\ \hline
 80\!&	0.60&	0.31\!&	0.98&	0.89&	0.99\!&	0.99&	1.00&	1.00\!\\ \hline
 90\!&	0.64&	0.36\!&	0.99&	0.93&	0.99\!&	0.99&	1.00&	1.00\!\\ \hline
100\!&	0.69&	0.41\!&	0.99&	0.95&	1.00\!&	0.99&   1.00&	1.00\!\\ \hline
\end{tabular}
\end{table}

\begin{table}
\caption{Empirical Power: Linear failure rate distribution}
\centering
\begin{tabular}{|l|l|l|l|l|l|l|l|l|}
\hline
\multirow{2}{*}{} & \multicolumn{2}{l|}{$\lambda=0.2$} & \multicolumn{2}{l|}{$\lambda=0.4$} & \multicolumn{2}{l|}{$\lambda=0.6$} & \multicolumn{2}{l|}{$\lambda=0.8$} \\ \cline{2-9}
 $n$& 5\% & 1\%& 5\% & 1\%& 5\%& 1\%&  5\% & 1\%  \\ \hline
 60\!&    0.49&	0.22\!&	0.68&	0.38&	0.79\!&	0.51&	0.87&	0.65\! \\ \hline
 70\!&	0.55&	0.27\!&	0.74&	0.46&	0.84\!&	0.61&	0.92&	0.74\! \\ \hline
 80\!&	0.60&	0.32\!&	0.80&	0.53&	0.89\!&	0.68&	0.94&	0.81\! \\ \hline
 90\!&	0.65&	0.36\!&	0.83&	0.59&	0.91\!&	0.74&	0.97&	0.86\! \\ \hline
100\!&	0.69&	0.41\!&	0.87&	0.65&	0.94\!&	0.80&	0.98&	0.90\! \\ \hline
\end{tabular}
\end{table}

\begin{table}
\caption{Empirical Power: Makeham distribution}
\centering
\begin{tabular}{|l|l|l|l|l|l|l|l|l|}
\hline
\multirow{2}{*}{} & \multicolumn{2}{l|}{$\lambda=0.2$} & \multicolumn{2}{l|}{$\lambda=0.4$} & \multicolumn{2}{l|}{$\lambda=0.6$} & \multicolumn{2}{l|}{$\lambda=0.8$} \\ \cline{2-9}
 $n$& 5\% & 1\%& 5\% & 1\%& 5\% & 1\%&  5\% & 1\%  \\ \hline

 60\! &	0.37 &	0.14 &	0.49\! &	0.22 &	0.65\! &	0.36 &	0.87  & 	0.63\! \\ \hline
70\! &	0.42 &	0.17 &	0.55\! &	0.26 &	0.72\! &  0.43 &	0.92 &	0.72\! \\ \hline
80\!&	0.46 & 	0.20 &	0.60\! &	0.31 &	0.78\! & 	0.49 &	0.94 &	0.79\! \\ \hline
90\!&	0.51& 0.23   &0.65\! &	0.35 &	0.82 &  0.56\! &	0.96 &	0.84\! \\ \hline
100\! &	0.55 &	0.27 &	0.70\! &	0.40 &	0.86\! &	0.62 &	0.98 &	0.90\! \\ \hline
\end{tabular}
\end{table}

\section{The case of censored observations}
Next we discuss how the censored observations can be incorporated in the proposed method.
Suppose we have randomly right-censored observations such that the censoring times are independent of the lifetimes. Under this set up the observed data are $n$ independent and identical copies of $(Y^*, \delta)$, with $Y^*=min(X^*,C)$, where $C$ is the censoring time and $\delta=I(X^*\leq C)$. Now we need to address the testing problem mentioned in Section 2 based on $n$ independent and identical observation $\{(Y_{i}^{*},\delta_i),1\leq i\leq n\}$. Observe that $\delta_i=1$  means $i^{th}$ object is not censored, whereas $\delta_i=0$ means that $i^{th}$  object is censored by $C$, on the right. Usually we need to redefine the measure $\Delta^*(F^*)$ to incorporates the censored observations.  The U-statistics formulations helps us to solve the problem in an easy way. Using the right-censored version of a U-statistic introduced by Datta et al. (2010)
 an estimator $\Delta(F^*)$ with censored observation is given by
\begin{equation*}\label{eq4.3}
\widehat{\Delta}_c=\frac{2}{n(n-1)}\sum_{i=1}^{n}\sum_{j<i;j=1}^{n}\frac{h(Y_{i}^*,Y_{j}^{*})\delta_i\delta_j}{\widehat{K}_{c}(Y_i^*)\widehat{K}_{c}(Y_j^{*})},
\end{equation*}
where $h(Y_{1}^*,Y_{2}^*)=\frac{1}{4}(4Y_{1}^*I(Y_{1}^*<Y_{2}^*)+4Y_{2}^*I(Y_{2}^*<Y_{1}^*)-Y_{1}^*-Y_{2}^*)$, provided $\widehat{K}_{c}(Y_i^*),\widehat{K}_{c}(Y_j^*)>0$, with probability 1 and $\widehat{K}_c$ is the Kaplan-Meier estimator of $K_c$, the survival function of the censoring variable C.
Similarly an estimator of $\mu^*$ is given by (Zhao and Tsiatis, 2000)
 \begin{equation*}
\widehat{X}_{c}^{*}=\frac{1}{n}\sum_{i=1}^{n}\frac{Y_i^* \delta_i}{\widehat{K}_{c}(Y_i^*)}.
\end{equation*}
Hence in right censoring situation, the test statistics is given by
\begin{equation*}
  \widehat{\Delta}_{c}^{*}=\frac{\widehat{\Delta}_c}{\widehat{X}_{c}^{*}},
\end{equation*} and  the test procedure is to reject $H_{0}$ in favour of $H_{1}$ for large values of $\widehat{\Delta}_c^*$.

Next we obtain the limiting distribution of the test statistic. Let $N_i^c(t)=I(Y_i^*\leq t, \delta_i=0)$ be the counting process corresponds to the censoring variable for the $i^{th}$ individual, $Z_i(t)=I(Y_i^*\geq t)$. Also let $\lambda_c$ be the  hazard rate of $C$.  The martingale associated with this counting process is given by
\begin{equation*}
M_i^c(t)=N_i^c(t)-\int_{0}^{t} Z_i(u) \lambda_c(u) du.
\end{equation*}
Let $G(x,y)=P(X_{1}^*\leq x, Y_{1}^*\leq y,  \delta=1), x\in \mathcal{X}$, $H(t)=P(Y_{1}^{*}\leq t)$ and
\begin{equation*}
w(t)=\frac{1}{\bar{H}(t)} \int_{\mathcal{X}\times[0,\infty)}{\frac{h_1(x)}{K_c(y-)}I(y>t)dG(x,y)},
\end{equation*}
where $h_1(x)=Eh(x,X_2^*).$ Next result follows from Datta et al. (2010) for the choice of the kernel $h(Y_{1}^*,Y_{2}^*)=\frac{1}{4}(4Y_{1}^*I(Y_{1}^*<Y_{2}^*)+4Y_{2}^*I(Y_{2}^*<Y_{1}^*)-Y_{1}^*-Y_{2}^*)$.
\begin{theorem}\label{thm5.1}
If $Eh^2(Y_1^*,Y_2^*)<\infty,$ $\int_{\mathcal{X}\times[0,\infty)}{\frac{h^{2}_1(x)}{K_c^2(y)}dG(x,y)}<\infty$ and  $\int_0^\infty w^2(t)\lambda_c(t)dt<\infty,$ then the distribution of $\sqrt{n}(\widehat{\Delta}_c-\Delta(F^*))$, as $ n \rightarrow \infty $, is Gaussian with mean zero and variance $4\sigma_{1c}^{2}$, where $\sigma_{1c}^2$  is given by
\begin{equation*}
\sigma_{1c}^{2}=Var\Big(\frac{h_1(X^*)\delta_1}{K_c(Y^{*}_1-)}+\int w(t) dM_1^c(t)\Big).
\end{equation*}
\end{theorem}
\begin{corollary}
Under the assumptions of Theorem \ref{thm5.1}, if $E(Y_1^2)<\infty,$   the  distribution of $\sqrt{n}(\widehat{\Delta}_c^*-\Delta^*(F^*))$, as $ n \rightarrow \infty $, is Gaussian with mean zero and variance $4\sigma_c^{2}$, where
\begin{equation}
\sigma_c^{2}=\frac{\sigma_{1c}^2}{\mu^{*2}}.
\end{equation}
\end{corollary}
\begin{corollary}Let $X^*$ be continuous non-negative random variable with $\bar{F^*}(x)=e^{-\frac{x}{\lambda} }$.  Under the assumptions of Theorem \ref{thm5.1}, if $E(Y_1^2)<\infty,$   the  distribution of $\sqrt{n}\widehat{\Delta}_c^*$, as $ n \rightarrow \infty $, is Gaussian with mean zero and variance $\sigma_{c0}^{2}$, where
\begin{equation}
\sigma_{c0}^{2}=\frac{4}{\lambda^{2}}Var\Big(\frac{(4F^*(X^*)-X^*)\delta_1}{4K_c(Y^{*}_1-)}+\int w(t) dM_1^c(t)\Big).
\end{equation}
\end{corollary}

Hence by Corollary 5.2, we know that the $\sqrt{n}\widehat{\Delta}_c^*$ has asymptotically  normal with mean zero and a variance that can be estimated by (5) and we denote it as $\widehat{\sigma}_{c0}$.  Hence we reject the null hypothesis in favour of $H_1$, if
\begin{equation*}
  \frac{\sqrt{n}\widehat{\Delta}_c^*}{\widehat{\sigma}_{c0}}\ge Z_{\alpha}.
\end{equation*}

Next we study the efficiency
loss due to censoring by computing the efficiency
of  our  test based on $\widehat{\Delta}^*$ for uncensored model
and the efficiency of the test based
on $\widehat{\Delta}^{*}_{c}$ for censored model. As both these tests have same asymptotic mean, the Pitman asymptotic relative efficiency (ARE) of the test
based on $\widehat{\Delta}^{*}_{c}$ with respect to the test based on $\widehat{\Delta}^*$  is given by
\begin{equation*}\label{eq6a}
  e=ARE(\widehat{\Delta}^{*}_{c},\widehat{\Delta}^*)= \frac{ \sigma^{2}_{0}}{\sigma_{c0}^{2}}.
\end{equation*}The quantity $(1-e)$ can  be taken as a measure of the efficiency loss (Lim and Park, 1993) due to censoring. From the above expression it is clear that the ARE value is independent of
the distributions belonging to the family of alternative hypothesis, but depends on the  distribution of $C$.
\begin{table}[h]
\caption{Asymptotic relative efficiency}
\centering
\begin{tabular}{|l|l|l|l|l|l|l|l|l|l|}
\hline
 $\lambda$& 0.5 & 0.4& 0.3& 0.2&  0.1& 0.05 &  0.01  \\ \hline
ARE\! &	 0.397\!& 0.433\!& 0.480\!& 0.547\!& 0.643\!& 0.700\!& 0.741\!
\\ \hline
\end{tabular}
\end{table}
Next we calculate the ARE value when the censoring variable $C$ has logistic distribution with distribution  function $ F(x)=\frac{1}{1+e^{-\frac{x}{\lambda }}}$. The ARE value for different values of $\lambda$  is given in Table 7. Table 7 shows that  as $\lambda$ decreases, the value of ARE increases and the efficiency loss decreases as the value of $\lambda$ (the amount of censoring) becomes small.


\end{document}